\documentclass[graybox]{svmult}

\usepackage{mathptmx}       
\usepackage{helvet}         
\usepackage{courier}        
\usepackage{type1cm}        
\usepackage{makeidx}         
\usepackage{graphicx}        
\usepackage{multicol}        
\usepackage[bottom]{footmisc}
\usepackage{amsmath, amsfonts, amssymb}
\usepackage{color}
\usepackage{tikz}

\usepackage{bm,enumitem}
\newcommand{\rhojp}{\rho_{j,+}}

\newcommand{\ujp}{u_{j,+}}

\newcommand{\pjp}{\tilde p_{j,+}}
\newcommand{\pjjm}{\tilde p_{j+1,-}}

\newcommand{\ps}{p_{j+1/2}^\star}
\newcommand{\us}{u_{j+1/2}^\star}
\newcommand{\usn}{u_{j+1/2}^{\star,n}}
\newcommand{\qs}{q_{j+1/2}^\star}
\newcommand{\ejp}{e_{j,+}}

\newcommand{\Ejp}{E_{j,+}}

\newcommand{\mjp}{m_{j,+}}
\newcommand{\hjp}{h_{j,+}}

\newcommand{\bu}{\bm{u}}

\newcommand{\bnu}{\bm{\nu}}
\newcommand{\bpi}{\bm{\pi}}

\makeindex             

\begin{document}

\title*{Lagrange-Flux schemes and the entropy property}
\titlerunning{Lagrange-flux schemes}
\author{Florian De Vuyst}

\institute{
Florian De Vuyst 
\at 
\'Ecole Normale Sup\'erieure Paris-Saclay \\
CMLA UMR 8536, 
61, avenue du Pr\'esident Wilson,
94235 Cachan cedex, France \\
\email{florian.de-vuyst@ens-paris-saclay.fr}}
%
\maketitle

\vspace{-1.0cm}
\abstract{
The Lagrange-Flux schemes are Eulerian finite volume schemes that make use of
an approximate Riemann solver in Lagrangian description with particular 
upwind convective fluxes. They have been recently designed as variant formulations
of Lagrange-remap schemes that provide better HPC performance on modern multicore
processors, see~[De Vuyst et al., OGST 71(6), 2016].
Actually Lagrange-Flux schemes show several advantages compared to Lagrange-remap schemes, especially for multidimensional problems: they do not require the computation of 
deformed Lagrangian cells or mesh intersections as in the remapping process.
The paper focuses on the entropy property of Lagrange-Flux schemes in their
semi-discrete in space form, for one-dimensional problems and for the compressible Euler equations as example. We provide pseudo-viscosity pressure terms
that ensure entropy production of order $O(|\Delta u|^3)$, where $|\Delta u|$ represents
a velocity jump at a cell interface. Pseudo-viscosity terms are also designed to vanish into  expansion regions as it is the case for rarefaction waves.
\keywords{Hyperbolic system, compressible Euler equations, finite volume, Lagrangian solver, Lagrange-remap scheme, Euler equations, discrete entropy property, numerical analysis
\\[5pt]
}
}
\vspace{-0.8cm}
\section{Introducing Lagrange-flux schemes}
%
%
Let us consider here the compressible Euler equations for two-dimensional problems. Denoting $\rho,\ \bu=(u_i)_i,$ $i\in\{1,2\}$,\ $p$ and $E$ the density,
velocity, pressure and specific total energy respectively, the mass, momentum and
energy conservation equations are
\begin{equation}
\partial_t U_\ell + \nabla\cdot(\bu\, U_\ell) + \nabla\cdot\bpi_\ell=0,
\quad \ell=1,\dots,4,
\end{equation}
where $U=(\rho,(\rho u_i)_i,\rho E)$, $\bpi_1=\vec 0$, $\bpi_2=(p,0)^T$,
$\bpi_3=(0,p)^T$ and $\bpi_4=p\bu$.
For the sake of simplicity, we will use a perfect gas equation of state
$p=(\gamma-1)\rho e$ with  $e=E-\frac{1}{2}|\bu|^2$ the internal energy
and $\gamma\in(1,3]$ the ratio of specific heats at constant volume. 
The speed of sound $c>0$ is such that $c^2=\frac{\partial p}{\partial \rho}|_s=\frac{\gamma p}{\rho}$.
The specific entropy~$s$ is given by $s=p/\rho^\gamma$. It is known that
the quantity $\eta(U):=-\rho\log(s)$ is a mathematical convex entropy for the system
and we look for physical weak entropy solutions $U$ that satisfy the inequality
\[
\partial_t\eta(U) + \partial_x(\eta(U) u) \leq 0
\]
in the sense of distributions \cite{GR96}. \medskip

Lagrange-Flux schemes have been derived in~\cite{FDV16,ECCOMAS2016,PONCET2016} from cell-centered Lagrange-{re}{map} schemes (see also~\cite{GASC}). Collocated Lagrangian solvers have been proposed a decade ago by Despr\'es-Mazeran~\cite{DESPRES} and by Maire et~al.~\cite{MAIRE}. By making the time 
step~$\Delta t$ tend to zero in Lagrange-remap schemes, it can be shown that this leads to the semi-discrete-in-space finite volume scheme (with standard notations):
\begin{equation}
\frac{d (U_\ell)_K}{dt} = -\frac{1}{|K|}\
\sum_{A\subset\partial K}|A|\, (U_\ell)_A^{upw}\, (\bu_A\cdot\nu_A)
- \frac{1}{|K|}
\sum_{A\subset\partial K}  |A|\, ((\bpi_\ell)_A\cdot \nu_A) 
\label{eq:11}
\end{equation}
for each $\ell=1,...,4$. In~\eqref{eq:11}, the notation $K$ stands for a generic control volume, $A$
is an edge of $K$, $\nu_A$ is the outward normal unit vector at the edge $A$, $(\bu_A\cdot\nu_A)$ is the normal fluid velocity at the edge~$A$ and $U_A^{upw}$ is the state at the
edge $A$, computed by some upwinding process.
%
We get a classical finite volume method in the form
\[
\frac{dU_K}{dt} = - \frac{1}{|K|}\,\sum_{A\subset \partial K} |A| \ \Phi_A
\]
with a numerical flux $\Phi_A$ whose components are
\begin{equation}
(\Phi_\ell)_A =  (U_\ell)_A^{upw}\, (\bu_A\cdot\nu_A) + (\bpi_\ell)_A\cdot \nu_A. 
\label{eq:12}
\end{equation}
Normal interface velocity $(\bu_A\cdot\bnu_A)$ and pressure $p_A$ can be computed by any approximate Riemann
solver in the Lagrangian frame (Lagrangian HLL solver~\cite{FDV16} for example) or derived
using a pseudo-viscosity approach. One can observe the simplicity
of expression~\eqref{eq:12} which is naturally consistent with the physical flux, and the way
pressure terms and convective terms are treated in a separate way.

Lagrange-Flux schemes have been since successfully extended to multi-material hydrodynamics
problems considering low-diffusive interface capturing schemes \cite{FDV2016b}.
\vspace{-0.8cm}
\section{The discrete entropy property}
%
In this section, we look for a discrete entropy property for particular Lagrange-Flux schemes
in their spatial semi-discrete form. For simplicity, we shall consider one-dimensional
formulations of the compressible Euler equations. The vector of conservative variables is $U=(\rho,\rho u,\rho E)$. We deal with conservative semi-discrete schemes in the form
\begin{equation}
\frac{dU_j}{dt} = -\frac{1}{|I_j|} \left( \Phi_{j+1/2} - \Phi_{j-1/2} \right)
\label{eq:fvm}
\end{equation}
for numerical fluxes in the form
\begin{equation}
\Phi_{j+1/2} =  U_{j+1/2}^{upw}\, \us + (0,\, \ps,\, \qs)^T.
\label{eq:phi}
\end{equation}
for each interval $I_j=(x_{j-1/2},x_{j+1/2})$, $x_{j+1/2}=(j+1/2)h$, where $h$ is the constant
space step, $\us$ and $\ps$ are interface velocity and pressure respectively, 
and~$\qs$ is consistent with the quantity $q=pu$ (not necessarily equal to~$\ps\us$).
Finally we will look for convected interface states $U_{j+1/2}^{upw}$ in the form
\[
U_{j+1/2}^{upw} = U_{j,+}\ 1_{(\us\geq 0)} + U_{j+1,-}\ 1_{(\us< 0)}.
\]
The interface quantities $\us$, $\ps$ and $\qs$ as well as $U_{j,+}$ and $U_{j+1,-}$  have to be determined in order to get the discrete entropy property.
The construction requires the analysis of a 1D cell-centered Lagrangian scheme first.
\vspace{-0.5cm}
\subsection{Lagrange step (semi-discrete version)}
%
We here follow ideas from Braeunig in the recent paper~\cite{JPB2016}. Each cell is splitted
up into two half-cells (denoted by $I_{j,+}$ and $I_{j,-}$ for the cell $j$) where the Euler equations are solved in each half-cell. Let us focus on the right half cell $I_{j,+}$
of~$I_j$.
Euler equations are discretized in the Lagrangian frame as follows: first,
the length of the half-cell $I_{j,+}$ evolves according to the geometric conservation law
$\frac{d\hjp}{dt} = 0$. The fluid mass $\mjp$ in the half-cell $I_{j,+}$ stays constant according to the mass conservation law: 
\[
\frac{d\mjp}{dt} = \frac{d(\hjp\rhojp)}{dt} = 0.
\]
Momentum and total energy conservation equations are written
\[
\mjp \frac{d\ujp}{dt} = -\left(\ps - p_j\right), \ \quad \mjp \frac{d\Ejp}{dt} = -\left(\qs - q_j\right)
\]
with $\Ejp = \frac{1}{2}(\ujp)^2+\ejp$. From the momentum equation, one can easily
get a balance equation for the kinetic energy
\[
\mjp \frac{d(\ujp)^2/2}{dt} = -\ujp\left(\ps - p_j\right)
\]
and then for the internal energy
\begin{equation}
\mjp \frac{d\ejp}{dt} + \left(\qs - q_j\right) - \ujp\left(\ps - p_j\right) = 0.
\label{eq:dte1}
\end{equation}
Of course, from the continuous partial differential equation $\rho D_t e + p \partial_x u = 0$, we would like
to reach a semi-discrete differential equation in the form
\begin{equation}
\mjp \frac{d\ejp}{dt} + \pjp \left(\us - u_j \right) = 0
\label{eq:dte2}
\end{equation}
for some pseudo-pressure $\pjp$ to be defined.
By identification between~\eqref{eq:dte1} and~\eqref{eq:dte2}, 
this leads to the following compatibility relation
\begin{equation}
\pjp \left(\us - u_j \right) + \ujp\left(\ps - p_j\right) = \qs - q_j.
\label{eq:c1}
\end{equation}
Similarly, into the left half-cell of cell $(j+1)$, we would get the compatibility
relation
\begin{equation}
\tilde p_{j+1,-} \left(u_{j+1}-\us\right) + u_{j+1,-}\left(p_{j+1}-\ps \right)
= q_{j+1} - \qs.
\label{eq:c2}
\end{equation}
At ``initial time'', it is natural to consider $u_{j,-} = u_{j,+} = u_j$ for
piecewise constant first-order discretizations. We will then consider this choice
in the sequel. So we recall the two compatibility relations subject to this choice:
\begin{eqnarray}
&& \pjp \left(\us - u_j \right) + u_j\left(\ps - p_j\right) = \qs - q_j, 
\label{eq:c1b}\\
&& \tilde p_{j+1,-} \left(u_{j+1}-\us\right) + u_{j+1}\left(p_{j+1}-\ps \right)= q_{j+1} - \qs.
\label{eq:c2b}
\end{eqnarray}
First, doing the difference between~\eqref{eq:c2b} and~\eqref{eq:c1b}, we get
an expression for $\qs$, depending on $\us$, $\ps$, $\pjp$ and $\pjjm$
(still to be determined).
Summing up~\eqref{eq:c1b} and~\eqref{eq:c2b} leads to a new compatibility relation
\begin{equation}
\boxed{
\pjp\left(\us - u_j\right)+\pjjm\left(u_{j+1}-\us\right) = \ps\left(u_{j+1}-u_j\right)}
\label{eq:c}
\end{equation}
that links the different local nonconservative products ``$p\Delta u$''. The question now
is to know how to define $\us$, $\ps$, $\pjp$ and $\pjjm$ in order to satisfy~\eqref{eq:c}
and achieve local entropy production in each half cell. There are (of course)
several candidates for that. The following result provides a simple choice for these
quantities, in the spirit of pseudo-viscosity terms.
\begin{theorem}\label{thm:1}
Consider the centered interface velocity 
\begin{equation}
\us = \frac{u_j+u_{j+1}}{2}
\label{eq:us}
\end{equation}
and the half cell pressure
\small
\begin{equation}
\pjp = \frac{p_j+\ps}{2} - \alpha\, (\rho c)_j\, \left(\us-u_j\right)_-
- \beta\, \rho_j\, |\us-u_j|\left(\us-u_j\right)_-
\label{eq:pjp}
\end{equation}
\normalsize
with the notation $x_-=\min(0,x)$, for some pseudo-viscosity constants $\alpha$ and $\beta$
(with a similar symmetric expression for $p_{j+1,-}$). Then
\begin{enumerate}
\item From the compatibility relation~\eqref{eq:c} we get the unique admissible value for 
the interface pressure~$\ps$:
\begin{equation}
\ps = \frac{\pjp+\pjjm}{2}
\label{eq:ps}
\end{equation}
\item We get the following half-cell entropy production
\begin{equation}
\mjp \frac{d\ejp}{dt} + \frac{p_j+\ps}{2}(\us-u_j) = \pi_{j,+}
\label{eq:tds}
\end{equation}
with
\small
\begin{equation}
\pi_{j,+} = \alpha\, (\rho c)_j\, [(\us-u_j)_-]^2 
+ \beta\, \rho_j \,|\us-u_j|\, [(\us-u_j)_-]^2 \geq 0.
\label{eq:pijp}
\end{equation}
\normalsize
For large velocity jumps $|u_{j+1}-u_j|$, we have
\[
\pi_{j,+} = O(|u_{j+1}-u_j|^3).
\]
\end{enumerate}
\end{theorem}
The proof of Theorem~1 is elementary and thus not detailed.
\begin{remark}
\begin{enumerate}[label=\alph*)]
\item Expressions~\eqref{eq:tds} and~\eqref{eq:pijp} show that $\dfrac{d s_{j,+}}{dt} \geq 0$
or equivalently
\begin{equation}
\frac{d}{dt} (-\mjp\log(s_{j,+})) \leq 0.
\label{eq:djlag}
\end{equation}
The order $O(|\Delta u|^3)$ for $\pi_{j,+}$ is the theoretically expected magnitude of entropy production. 
\item Remark that for expansion regions like rarefaction fans, we have $\pi_{j,+}=0$,
meaning that there is not artificial viscosity in this case.
\item Dimensionless pseudo-viscosity constants $\alpha$ and $\beta$ are free parameters
let to the user. If we try to connect $\ps$ to a mean pressure returned by some approximate Riemann solvers (like HLL Lagrange or acoustic solver, see~\cite{Toro09}), then one finds that 
$\alpha=\frac{1}{2}$ and $\beta$ should be of order $1$.
\item Of course, because we have entropy production into each half cell $I_{j,+}$
and $I_{j,-}$, we have also entropy production in the whole cell $I_j$.
\item Because of the two distinct evolutions on each half-cells, a full space-time explicit discretization would involve a CFL condition less than $\frac{1}{2}$ as a necessary
stability condition.
\end{enumerate}
\end{remark}
\vspace{-0.5cm}
\subsection{Remap step}
%
If the Lagrangian scheme above is fully discretized, of course after a time step $\Delta t^n$
we have a discrete updated solution $U^{n+1,L}(.)$ defined on a deformed Lagrangian mesh
(the superscript $L$ stands for ``Lagrange''). 
A remapping process is needed to remap the discrete solution on the reference Eulerian (fixed) mesh.

The remapping step is nothing else but a conservative projection of the conservative quantities defined in the deformed Lagrangian mesh onto the reference (Eulerian) mesh. For a first-order accurate remapping procedure, we have only to perform averages
on piecewise constant functions. For any convex function $\eta$, using Jensen's inequality,
we have
\begin{equation}
\eta(U_j^{n+1})=
\eta\left(\frac{1}{|I_j|}\int_{I_j} U^{n+1,L}(x)\, dx\right) \ 
\leq \ \frac{1}{|I_j|}\int_{I_j} \eta(U^{n+1,L}(x))\, dx.
\label{eq:jensen}
\end{equation}
As soon as the time step $\Delta t^n$ is chosen in order not to create unexpected inverted cells (according to some ``material'' CFL-like condition), formula~\eqref{eq:jensen} also leads to the following entropy inequality in conservative form
\begin{equation}
\eta(U_j^{n+1}) \leq \eta(U_j^{n+1,L}) 
- \frac{\Delta t^n}{|I_j|} \left(\Psi_{j+1/2}^{n+1,L}-\Psi_{j-1/2}^{n+1,L}\right)
\label{eq:ineq}
\end{equation}
with the numerical entropy flux
\[
\Psi_{j+1/2}^{n+1,L} =  \eta(U_{j,+}^{n+1,L})\,(\usn)_+ 
                     + \eta(U_{j+1,-}^{n+1,L})\,(\usn)_-
\]
which is consistent with the physical entropy flux $\Psi=\eta u$. 
%
\vspace{-0.5cm}
\subsection{Back to the semi-discrete Eulerian Lagrange-Flux scheme}
%
Following the construction of the Lagrange-flux schemes presented in~\cite{FDV16},
we now consider the Lagrange-remap scheme described in the two above subsections 
and observe what happens when $\Delta t^n$ tends to 0. We get a semi-discrete
Lagrange-flux scheme written in the form~\eqref{eq:fvm},\eqref{eq:phi}.
From~\eqref{eq:ineq}, at the limit $\Delta t^n\rightarrow 0$ we get the
semi-discrete entropy inequality
\begin{eqnarray*}
\frac{d}{dt}\eta(U_j) + \frac{1}{|I_j|} \left(\Psi_{j+1/2}-\Psi_{j-1/2}\right) \leq 0
\label{eq:ineq2}
\end{eqnarray*}
with the numerical entropy flux
\begin{equation}
\Psi_{j+1/2} =  \eta(U_{j,+})\,(\us)_+  \,+\, \eta(U_{j+1,-})\, (\us)_-.
\label{eq:21}
\end{equation}
\vspace{-1.2cm}
\subsection{Discussion}
%
\begin{enumerate}
\item The analysis calls for the use of entropy balances in half-cells for
achieving a discrete entropy property. This is quite similar to the use of subcell
evolutions in the framework of cell-centered Lagrangian schemes proposed in the GLACE scheme by Despr\'es-Mazeran~\cite{DESPRES} or the EUCCLHYD scheme by Maire et al.~\cite{MAIRE}. 
\item The entropy result in this paper is obtained for a Lagrange-flux scheme only in its semi-discretized form. We are aware that, for a full-discretized entropy inequality result, we have to pay attention about time discretization. In our opinion, some variables have to be treated in an implicit way to achieve a full discrete entropy inequality. 
\item Other choices than~\eqref{eq:us}, \eqref{eq:pjp} and~\eqref{eq:ps} can be derived
to get~\eqref{eq:c}, but are not discussed in the present paper. 
\item The choice of interface velocity~\eqref{eq:us} is somewhat surprising because
it is a centered discretization, and does not connect to the interface velocity
returned by an approximate Riemann solver like HLL-Lagrange or the Lagrangian acoustic 
solver. Of course, one can add an additional artificial viscosity term into~\eqref{eq:us}
for even-better stability, but we may lose the theoretical  property that 
$\pi_{j,+}\geq 0$. 
\item It is clear that, for a full discretized Lagrange-flux scheme, the positivity of the density variable will be preserved during time iterations under suitable CFL condition. The issue of the positivity preservation for the internal energy will be discussed in a future paper.
\end{enumerate}
\vspace{-0.8cm}
\section{Numerical experiments}
\begin{figure}[h!]
\hspace{-1.4cm}\includegraphics[width=1.2\textwidth]{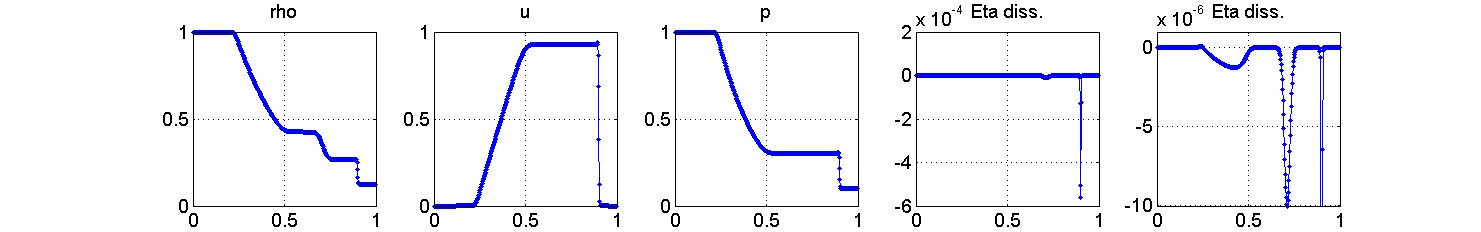}
\caption{Sod's 1D shock tube problem discrete solution at time $T=0.23$.
Number of cells: 400.}\label{fig:1}
\end{figure}

\begin{figure}[h!]
\hspace{-1.2cm}\includegraphics[width=1.2\textwidth]{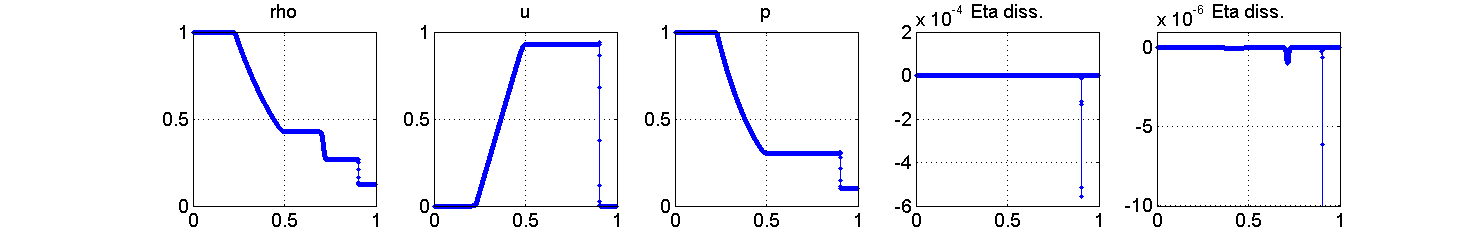}
\caption{Sod's 1D shock tube problem discrete solution at time $T=0.23$.
Number of cells: 4000.}\label{fig:2}
\end{figure}

%
As an example, we show how the Lagrange-flux scheme in its full discretized,
first order and pure explicit
form behaves on the reference one-dimensional Sod's shock tube problem of left state 
$(\rho,u,p)_L=(1,0,1)$ and right state $(\rho,u,p)_R=(0.125,0,0.1)$.
For pseudo-viscosity coefficients, we use $\alpha=\frac{1}{2}$ and $\beta=\frac{\gamma+1}{2}$.
A uniform grid mesh made of $N$ cells is used with step size $h=\frac{1}{N}$. Computations are performed under the CFL condition $\frac{1}{4}$. We also compute the entropy production defined as
\[
\Pi_j^{n,n+1} = \eta(U_j^{n+1})-\eta(U_j^n) + \frac{\Delta t^n}{h}\left(\Psi_{j+1/2}^n-\Psi_{j-1/2}^n\right)
\]
with the numerical entropy flux~\eqref{eq:21}. For Figure~\ref{fig:1} we use $N=400$ mesh
cells and for Figure~\ref{fig:2} we use $N=4000$ cells. For each figure, we show the discrete
solution at final time $T=0.23$ including density, velocity, pressure and entropy dissipation profiles. During the simulation, one can observe a very slight entropy dissipation defect at the
top of the rarefaction fan, otherwise $\Pi_j^{n,n+1}$ has the right sign elsewhere.
\vspace{-0.5cm}
\vspace{-0.5cm}
\vspace{-1.3cm}
\section{Concluding remarks and perspectives}
%
\vspace{-0.3cm}
In this paper, we have justified the new family of Lagrange-flux schemes
from the point of view of the discrete entropy stability property. We have found
pseudo-viscosity terms that return an entropy production at the right order of
magnitude for a semi-discretized version of Lagrange-Flux schemes.
Even if further analysis is needed (multidimensional extension, analysis of the full space-time discretized scheme), we believe that Lagrange-flux schemes are very promising because of their
generic formulation, simplicity of implementation, simple expressions of numerical fluxes, flexibility for achieving higher-order accuracy, HPC performance (see the references below),
and flexibility for multiphysics coupling. 
%
%
\vspace{-0.4cm}
\begin{acknowledgement}
The work is part of the \emph{LRC MESO} joint lab between CEA DAM DIF and CMLA. The author would like to thank Dr. Jean-Philippe Braeunig for valuable discussions on this subject.
The author also thanks the anonymous reviewers for their constructive comments.
\end{acknowledgement}
%
%


\vspace{-0.9cm}
\bibliographystyle{spmpsci}
\bibliography{biblio}

\end{document}